\documentstyle[11pt]{article}

\def\ni{\noindent}

\def\b{\beta}

\begin{document}

\rightline{ITP-97-8E \hspace{8mm}}
\rightline{q-alg/9709036 \hspace{3mm}}
\vspace{28mm}

\centerline
{\bf $q$-DEFORMED ALGEBRAS $U_q({\rm so}_n)$ AND
                                 THEIR REPRESENTATIONS  }
\vspace{18mm}

\centerline
{ A. M.~GAVRILIK\footnote{ {\it E-mail address:} omgavr@gluk.apc.org}
                         and N. Z.~IORGOV }

\vspace{12mm}
\centerline{
Institute for Theoretical Physics of National Academy of
Sciences of Ukraine,}
\centerline
{Metrologichna 14b, 252143 Kiev, Ukraine}

\vspace{30mm}

\begin{abstract}
   For the nonstandard $q$-deformed algebras $U_q({\rm so}_n),$
defined recently in terms of trilinear relations for generating
elements, most general finite dimensional irreducible
representations directly corresponding to those of nondeformed
algebras ${\rm so}(n)$ (i.e., characterized by the same
sets of only integers or only half-integers as in highest weights
of the latter) are given explicitly in a $q$-analogue of
Gel'fand-Tsetlin basis. Detailed proof, for $q$ not equal to
a root of unity, that representation operators indeed satisfy
relevant (trilinear) relations and define finite dimensional
irreducible representations is presented.  The results show
perfect suitability of the Gel'fand-Tsetlin formalism
concerning (nonstandard) $q$-deformation of  ${\rm so}(n)$ .
\end{abstract}

\vfill

\newpage
\centerline {\sc 1. Introduction}
\medskip

The so called {\it standard} quantum deformations
$U_q({\rm B}_r)$,\ $U_q({\rm D}_r)$
of Lie algebras of the orthogonal groups ${\rm SO}(2r+1)$
,\  ${\rm SO}(2 r)$
were given by Drinfeld [1] and Jimbo [2]. However, there are some
 reasons why to
consider a ${q}$-deformation of the orthogonal Lie algebras
 ${\rm so}(n,{\bf C})$
other than the standard quantum algebras
$U_q({\rm B}_r)$,\ $U_q({\rm D}_r)$.
In order to construct explicitly (possibly, all) finite dimensional
irreducible representations of
$q$-deformed algebras $U_q({\rm so}(n,{\bf C}))$ and their
 compact real forms
$U_q({\rm so}_n)$, one needs a $q$-analogue of the
Gel'fand--Tsetlin (GT)
basis and the GT action formulas [3],[4] which requires existence of
canonical embeddings $q$-analogous to the chain ${\rm so}(n,{\bf C})
 \supset
{\rm so}(n-1,{\bf C})\supset \cdots \supset {\rm so}(3,{\bf C})$.
 Obviously, such a chain
does not exist for the standard quantum algebras $U_q({\rm B}_r)$
and $U_q({\rm D}_r)$, only single-type embedding
$U_q({\rm B}_{r-1}) \subset U_q({\rm B}_r)$ and embedding
$U_q({\rm D}_{r-1}) \subset U_q({\rm D}_r)$ do hold. Moreover,
 as it was pointed
out by Faddeev, Reshetikhin and
Takhtadzyan [5], standard quantum orthogonal groups
${\rm SO}_q(n,{\bf C})$
defined on the base of $R$-matrix do not possess, when $n>3$, real
forms of Lorentz signature (which are necessary for investigating
important aspects of quantized $n$-dimensional space-times).
Again, in developing representation theory of a $q$-analogue of the
Lorentz algebras ${\rm so}(n,1)$ in a way parallel to the classical
 $(q=1)$
case, see [6] and references therein, one also needs $q$-analogues of
the GT basis and the GT formulas which certainly require that
the chain of inclusions $U_q({\rm so}(n,1)) \supset U_q({\rm so}(n))
\supset U_q({\rm so}(n-1)) \supset \cdots \supset U_q({\rm so}(3))$
should exist.

An alternative approach to $q$-deformation of the
orthogonal and pseudoorthogonal Lie algebras was proposed in
[7],[8] and further developed in [9]-[11].
Such {\it nonstandard} $q$-analogue $U_q({\rm so}(n,{\bf C}))$ of the
 Lie  algebra
${\rm so}(n,{\bf C})$ is constructed not in terms of Chevalley basis
(i.e., Cartan subalgebra and simple root elements, as it was done
within standard
approach of [1],[2]) but, starting from ${\rm so}(n,{\bf C})$
formulated as the complex associative algebra with ${n-1}$
generating elements $I_{21}, I_{32}, \ldots , I_{n,n-1}$.  The {\it
nonstandard} $q$-analogues $U_q({\rm so}(n,{\bf C}))$ are formulated in
a uniform fashion for all values $n\geq 3$, and guarantee validity of
the canonical chain of embeddings $$ U_q({\rm so}(n,{\bf C})) \supset
U_q({\rm so}(n-1,{\bf C}))\supset \cdots
 \supset U_q({\rm so}(3,{\bf C})).          \leqno(1)
$$
Due to this, one can attempt to develop a $q$-analogue
of GT formalism. Other viable feature is the
admittance by these $U_q({\rm so}(n,{\bf C}))$ of all the noncompact
real forms corresponding to those in classical case. In particular,
validity of inclusions $U_q({\rm so}(n,1)) \supset U_q({\rm so}(n))
\supset U_q({\rm so}(n-1)) \supset ... \ \ $ can be exploited
in order to analyze infinite dimensional representations of
$q$-deformed Lorentz algebras (see [10] where detailed study of
class 1 representations of
$U_q({\rm so}(n,1))$ is presented).

The purpose of this paper is to provide,
{\it along with detailed proof of validity}, exact matrix
realization of finite dimensional irreducible representations
of the nonstandard $q$-deformation $U_q({\rm so}(n,{\bf C}))$
 that directly correspond to the finite dimensional
irreducible representations
of ${\rm so}(n,{\bf C})$ given by highest weights with all integral or
all half-integral components. Those representations, with an additional
restriction on the deformation parameter $q,$ yield irreducible
infinitesimally unitary (or $*$-) representations of the 'compact'
real form $U_q({\rm so}_{n})$.

\bigskip
\centerline
{\sc  2. Nonstandard $q$-deformed algebras $U_q({\rm so}_{n})$}
\bigskip

According to [7]-[10], the nonstandard
$q$-deformation $U_q({\rm so}(n,{\bf C}))$ of the Lie  algebra
 ${\rm so}(n,{\bf C})$ is
given as a complex associative algebra with ${n-1}$
generating elements $I_{21}, I_{32}, \ldots ,$ $I_{n,n-1}$
obeying the defining relations (denote $q+q^{-1}\equiv [2]_q$)
\newpage
$$
I_{k,k-1}^2I_{k-1,k-2} + I_{k-1,k-2}I_{k,k-1}^2 -
[2]_q \ I_{k,k-1}I_{k-1,k-2}I_{k,k-1} = -I_{k-1,k-2}, \leqno(2)
$$
$$
I_{k-1,k-2}^2I_{k,k-1} + I_{k,k-1}I_{k-1,k-2}^2 -
[2]_q \ I_{k-1,k-2}I_{k,k-1}I_{k-1,k-2} = -I_{k,k-1},    \leqno(3)
$$
$$
[I_{i,i-1},I_{k,k-1}] =0  \qquad {\rm if} \quad \mid {i-k}\mid >1,
                                                     \leqno(4)
$$
where $i,k=2,3,...,n$. These
relations are certainly true at $q = 1$, i.e., within ${\rm so}(n,{\bf C})$.
The real forms - compact $U_q({\rm so}_n)$ and noncompact
$U_q({\rm so}_{n-1,1})$ -
are singled out from the algebra $U_q({\rm so}(n,{\bf C}))$ by imposing the
$*$-structures
$$I_{k,k-1}^* = - I_{k,k-1}, \quad k=2,...,n, \leqno(5)
$$
and
$$
I_{k,k-1}^*=-I_{k,k-1}, \quad k=2,...,n-1, \qquad I_{n,n-1}^*=I_{n,n-1},
                                                    \leqno(6)
$$
respectively\footnote {Other noncompact real forms
$U_q({\rm so}_{n-p,p}),$
defined by their corresponding $*$-structures,
are also possible [7],[10].}.
It is obvious that complex as well as compact and noncompact
real $q$-deformed algebras so defined are compatible with the canonical
reductions mentioned above. This enables one to apply a $q$-analogue
of GT formalism in constructing representations. Finite dimensional
representations $T$ of the algebras $U_q({\rm so}_n)$ are
characterized by signatures completely analogous to highest weights of
the algebras ${\rm so}(n)$ and described by means of action formulas
for the operators $T(I_{k,k-1})$, $k=2,...,n$, which satisfy the
relations (1)-(3) and the relations $$T(I_{k,k-1})^* = - T(I_{k,k-1}).
\leqno(7)
$$

Finite dimensional representations of $U_q({\rm so}_{3}),$
$U_q({\rm so}_{4})$ and infinite dimensional representations of their
noncompact analogues $U_q({\rm so}_{2,1}),$ $U_q({\rm so}_{3,1})$ were
studied in [9],[10].

\medskip
{\it Remark. } As pointed out in [9]-[10], the algebra
$U_q({{\rm so}}_n)$
for $n=3,$ when presented in terms of {\it bilinear} defining
relations (i.e., with $q$-deformed commutators) for
 $I_{21}$ and $I_{32}$ is isomorphic to
the (cyclically symmetric, Cartesian) $q$-deformed algebra which was
studied by Odesskii [13] and Fairlie [14] where some representations
 were also given. Classification of irreducible $*$-representations of
 the algebra $U_q({\rm so}_3)$ is worked out in [15].

\bigskip

\centerline
{\sc 3. Representations of the algebras $U_q({{\rm so}}_n)$}
\medskip

In this section we describe explicitly finite dimensional
representations of the algebras $U_q({\rm so}_{n}), \ n \ge 3$
in the framework of a $q$-analogue of GT formalism.
These are given by signatures - sets ${\bf m}_{n}$
consisting of
$[\frac{n}{2}]$ components $m_{1,n}, m_{2,n},..., m_{\left
[\frac{n}{2}\right ],n}$ (here $[\frac{n}2]$ denotes integer part
of $\frac{n}2$) that satisfy the dominance condition,
respectively, for $n=2p+1$ and $n=2p$ :
$$
m_{1,2p+1}\ge m_{2,2p+1}\ge ... \ge
m_{p,2p+1}\ge 0 ,                                  \leqno(8a)
$$
$$
m_{1,2p}\ge m_{2,2p}\ge ... \ge m_{p-1,2p}\ge |m_{p,2p}|.
                                                    \leqno(8b)
$$
For a basis in representation space
we exploit the $q$-analogue of GT basis [3],[4].  Its elements are
labelled by GT schemes
$$
\left.  \{\xi_{n} \}
\equiv \left\{ \begin{array}{l} {\bf m}_{n} \\ {\bf m}_{n-1} \\ \dots
\\ {\bf m}_{2} \\ \end{array}
\right.  \right\}
\equiv \{ {\bf m}_{n},\xi_{n-1}\}\equiv \{{\bf m}_{n} ,
{\bf m}_{n-1} ,\xi_{n-2}\} ,
                                                           \leqno(9)
$$
 where the components of ${\bf m}_{n}$ and ${\bf m}_{n-1}$ satisfy the
"betweenness" conditions
$$
m_{1,2p+1}\ge m_{1,2p}\ge m_{2,2p+1} \ge m_{2,2p} \ge ...
\ge m_{p,2p+1} \ge m_{p,2p} \ge -m_{p,2p+1}  ,
                                        \leqno(10a)
$$
$$
m_{1,2p}\ge m_{1,2p-1}\ge m_{2,2p} \ge m_{2,2p-1} \ge ...
\ge m_{p-1,2p-1} \ge \vert m_{p,2p} \vert .
                                        \leqno(10b)
$$
Basis element defined by scheme $\{\xi_{n} \}$ is denoted
as $\vert \{\xi_{n} \} \rangle $ or simply as $\vert
\xi_{n} \rangle $.

We use standard denotion for $q$-number
$[x]\equiv (q^x-q^{-x})/(q-q^{-1})$ corresponding to an integer
or a half-integer $x$.
Also, it is convenient to introduce
the so-called $l$-coordinates
$$
l_{j,2p+1}=m_{j,2p+1}+p-j+1,  \hspace{12mm}
                                  l_{j,2p}=m_{j,2p}+p-j .   \leqno(11)
$$
Infinitesimal operator $I_{2p+1,2p}$ of the representation, given by
${\bf m}_{2p+1}$,
of $U_q({\rm so}_{2p+1})$ acts upon GT basis elements,
labelled by (9), according to (here $\b\equiv\xi_{2p-1}$)
$$
I_{2p+1,2p}\vert {\bf m}_{2p+1}, {\bf m}_{2p},\b\rangle=
\sum^p_{j=1}A^j_{2p}({\bf m}_{2p})
            \vert {\bf m}_{2p+1},{\bf m}^{+j}_{2p},\b\rangle
                                        \hskip 32mm$$
$$
\hspace{10mm} -\sum^p_{j=1}A^j_{2p}({\bf m}^{-j}_{2p})
\vert {\bf m}_{2p+1},{\bf m}^{-j}_{2p},\b\rangle         \leqno(12)
$$
and the operator $I_{2p,2p-1}$ of the representation, given by
${\bf m}_{2p}$, of $U_q({\rm so}_{2p})$ acts as
(here $\b\equiv\xi_{2p-2}$)
$$
I_{2p,2p-1}\vert {\bf m}_{2p}, {\bf m}_{2p-1},\b\rangle=
\sum^{p-1}_{j=1}B^j_{2p-1}({\bf m}_{2p-1})
\vert {\bf m}_{2p},{\bf m}^{+j}_{2p-1},\b\rangle
                                        \hskip 29mm
$$
$$                  \hspace{10mm}
-\sum^{p-1}_{j=1}B^j_{2p-1}({\bf m}^{-j}_{2p-1})
\vert {\bf m}_{2p},{\bf m}^{-j}_{2p-1},\b\rangle
$$
$$                  \hspace{11mm}
+ \ i \ C_{2p-1}({\bf m}_{2p-1})
\vert {\bf m}_{2p},{\bf m}_{2p-1},\b\rangle .              \leqno(13)
$$
In these formulas, ${\bf m}^{\pm j}_n$ means that the $j$-th component
$m_{j,n}$ in signature ${\bf m}_n$ is to be replaced
by $m_{j,n}\pm 1$; matrix elements
$A^j_{2p}, \ $ $B^j_{2p-1}, \ $ $C_{2p-1}$ are given by the expressions
\medskip
$$
A^j_{2p}(\xi_{2p+1}) =    d(l_{j,2p})
\left|
\frac{ \prod_{i=1}^p
[l_{i,2p+1}+l_{j,2p}][l_{i,2p+1}-l_{j,2p}-1] }
 {  \prod_{i\ne j}^p
[l_{i,2p}+l_{j,2p}] [l_{i,2p}-l_{j,2p}] }  \right.
                       \hspace{28mm}
$$
$$  \hspace{39mm} \times \left. \frac{   \prod_{i=1}^{p-1}
[l_{i,2p-1}+l_{j,2p}][l_{i,2p-1}-l_{j,2p}-1] }
 {  \prod_{i\ne j}^p
[l_{i,2p}+l_{j,2p}+1] [l_{i,2p}-l_{j,2p}-1]  }
\right|^{\frac12}                                       \leqno(14)
$$
\vspace{1mm}
with
$$
  d(l_{j,2p})\equiv  \left( {[l_{j,2p}] [l_{j,2p}+1]
\over [2 l_{j,2p}] [2 l_{j,2p}+2]} \right)^{\frac12}     \leqno(15)
$$
and
$$
B^j_{2p-1}(\xi_{2p})=       \left|
\frac { \prod_{i=1}^p[l_{i,2p}+l_{j,2p-1}][l_{i,2p}-l_{j,2p-1}] }
{   [2 l_{j,2p-1}+1] [2 l_{j,2p-1}-1]
   \prod_{i\ne j}^{p-1}
 [l_{i,2p-1}+l_{j,2p-1}] [l_{i,2p-1}-l_{j,2p-1}] }
\right.  \
$$
$$   \hspace{18mm}   \times \left.
\frac { \prod_{i=1}^{p-1}
[l_{i,2p-2}+l_{j,2p-1}][l_{i,2p-2}-l_{j,2p-1}] }
  { [l_{j,2p-1}]^2   \prod_{i\ne j}^{p-1}
[l_{i,2p-1}+l_{j,2p-1}-1] [l_{i,2p-1}-l_{j,2p-1} -1]  }
\right|^{\frac12} ,                                      \leqno(16)
$$
\vspace{1mm}
$$
C_{2p-1}(\xi_{2p}) = \frac{
\prod_{s=1}^p [ l_{s,2p} ]
\prod_{s=1}^{p-1} [ l_{s,2p-2} ] }
  { \prod_{s=1}^{p-1} [l_{s,2p-1}] [l_{s,2p-1} - 1] } .    \leqno(17)
$$
From (17) it is seen that $C_{2p-1}$ in (13) identically
vanishes if $m_{p,2p}=l_{p,2p}=0.$

\smallskip
\ni{\it Remark.}
Matrix elements $B^j_{2p-1}(\xi )$ and $C_{2p-1}(\xi )$ are
nothing but "minimal" deformation of their classical (i.e., $q=1$)
counterparts. On the other hand, because of appearance of the nontrivial
factor $d(l_{j,2p})$ (which replaces the constant multiplier $\frac12$ of
classical case) in $A^j_{2p}(\xi )$, these matrix elements deviate from
"minimal" deformation.

\ni {\bf Proposition.}
{\it Let $q^N \ne 1,\ N\in {\bf Z}\backslash \{0\}$, $q \in {\bf C}$.
Representation operators $T_{{\bf m}_n}(I_{k,k-1}),$
$\ k=2,...,n,$ of
the representation $T_{{\bf m}_n},$ characterized by the signature
${\bf m}_n$ with all integral or all half-integral components
satisfying (8),
given in the $q$-analogue of GT basis (9),(10) by the action
formulas (12)-(17) satisfy the defining relations (2)-(4) of the
algebra} $U_q({\rm so}(n,{\bf C}))$ {\it for both even and odd $n$
and define finite dimensional irreducible representations
of this algebra.
If in addition $q=e^h$ or $q=e^{{\rm i}h},$ $\ h \in {\bf R},$ these
operators satisfy the condition (7) and thus provide
(infinitesimally unitary or) $*$-representations of
$U_q({\rm so}_n)$.}

\bigskip
\ni {\it Proof. }
First, we prove that the representation operators
$T_{{\bf m}_n}(I_{k,k-1}),\ \ k=2, \ldots ,n $ of the algebra
$U_q({\rm so}(n,{\bf C}))$
given in the $q$-analog of GT basis by explicit formulas (12)-(17)
indeed satisfy the defining relations (2)-(4). The fact that any pair
of operators $T_{{\bf m}_n}(I_{k,k-1})$ and $T_{{\bf
m}_n}(I_{i,i-1})$,\footnote { For convenience, below by the
symbol $I_{k,k-1}$ we mean also the corresponding representation
operator $T_{{\bf m}_n}(I_{k,k-1})$.} where $\vert i-k \vert > 1$ do
commute (cf.(4)) follows from the structure of their action formulas (
e.g. in the "nearest" case of $i-k=2$ the operator $I_{k,k-1}$ changes
components of signatures ${\bf l}_{k-1}$ while its matrix
elements do not depend on ${\bf l}_{k+1}$  affected by the action of
$I_{k+2,k+1}$, and vice versa). Now let us prove the relations
(2)-(3) for $I_{k,k-1}$ taken for sequential (odd and even)
values of $k$, say $2p+1$ and $2p$:
\smallskip
$${\bf I.} \ \               \hskip 9mm
I_{2p+1,2p}^2 I_{2p,2p-1} +  I_{2p,2p-1}I_{2p+1,2p}^2 - [2]\
I_{2p+1,2p} I_{2p,2p-1} I_{2p+1,2p}     \hspace{7mm}
$$
$$= -
I_{2p,2p-1};$$
$${\bf II.}\ \               \hskip 8mm
I_{2p,2p-1}^2 I_{2p+1,2p} +  I_{2p+1,2p}I_{2p,2p-1}^2 -
[2]\ I_{2p,2p-1} I_{2p+1,2p} I_{2p,2p-1}      \hspace{8mm}
$$
$$= - I_{2p+1,2p};
$$
$${\bf III.} \ \             \hskip 4mm
I_{2p,2p-1}^2 I_{2p-1,2p-2} +  I_{2p-1,2p-2}I_{2p,2p-1}^2 -
[2]\ I_{2p,2p-1} I_{2p-1,2p-2} I_{2p,2p-1}
$$
$$
 = - I_{2p-1,2p-2};
$$
$${\bf IV.}\ \
I_{2p-1,2p-2}^2 I_{2p,2p-1} +  I_{2p,2p-1}I_{2p-1,2p-2}^2 -
[2]\ I_{2p-1,2p-2} I_{2p,2p-1} I_{2p-1,2p-2}
$$
$$ = - I_{2p,2p-1}.
$$

Action of left and right hand sides of each of these equalities upon
some generic basis vector $\mid \xi \rangle$ of the representation
space ${\cal V}_{{\bf m}_n}$ produces, besides the initial basis vector
(if any), also a number of other basis vectors. It is necessary to
examine all the encountered resulting basis vectors for each of
the relations (I)-(IV).

\centerline
{\it Case of (I)}

\smallskip
The list of resulting basis vectors for (I) is as follows.

\ni Action of RHS of (I) on the vector  $\mid \xi \rangle$
leads to vectors :

\smallskip
(I.1.a)  $\mid \xi \rangle$ (unchanged vector):  \hspace{46mm}
1 vector

\smallskip
(I.1.b) $\vert l_{j,2p-1}\pm 1\rangle $:      \hspace{52mm}
$2(p-1)$  vectors

\smallskip

\ni Action of LHS of (I) on the same vector  $\mid \xi \rangle$
leads to vectors :

\smallskip
(I.1.a)  $\mid \xi \rangle$ (unchanged vector):  \hspace{46mm}
1 vector

\smallskip
(I.1.b) $\vert l_{j,2p-1}\pm 1\rangle $:   \hspace{52mm}
$2(p-1)$  vectors

\smallskip
(I.2.a) $\vert l_{j',2p}\pm 2\rangle$:    \hspace{64mm} $2p$  vectors

\smallskip
(I.2.b) $\vert l_{j',2p}\pm 2;l_{j,2p-1}\pm 1\rangle$:
                                       \hspace{36mm}   $4p(p-1)$ vectors

\smallskip
(I.3.a) $\vert l_{j',2p}\pm 1;l_{j'',2p}\pm 1\rangle$:
                                     \hspace{38mm} $2p(p-1)$  vectors

\smallskip
(I.3.b) $\vert l_{j',2p}\pm 1;l_{j'',2p}\pm 1;l_{j,2p-1}\pm 1\rangle$:
                                \hspace{12mm}   $4p(p-1)(p-1)$  vectors

\smallskip
\ni(only the changed component(s) off those labelling the vectors
are indicated; the signs $\pm$ in the last three items take
their values independently).
\medskip

Comparison of expressions at the same fixed vector $\mid \tilde\xi \rangle$
in the LHS and RHS yields some particular relation (in what follows,
we use for it the term {\it relation corresponding to the vector}
$\mid \tilde\xi \rangle$). The proof of  equality (I) will be
achieved if all the possible relations (obtained from examination
of the complete list of the encountered basis vectors for this
equality) are proven to hold identically. For every type of the vectors
in the list given above, we prove validity of relation corresponding to
one typical representative; the proof of relations corresponding to
the other vectors of this same type is completely analogous. For
instance, concerning the type (I.3.a) we prove the relation
corresponding to typical representative
$\vert l_{j',2p}+1;l_{j'',2p}+1\rangle$.

The vectors (I.2.a), (I.2.b), (I.3.a), (I.3.b) appear only in the LHS.
Hence, the relations which correspond to them have zero right hand sides;
verification of these relations proceeds in straightforward way.
Let us start with these 4 cases.

\smallskip
\centerline
{\it Relations corresponding to (I.2.a)}
\smallskip

   To typical vector $\vert l_{j',2p}+2 \rangle $ there corresponds
the relation
$$
C_{2p-1}(l_{j',2p}+2)+C_{2p-1}(l_{j',2p})-[2] C_{2p-1}(l_{j',2p}+1)=0.
$$
On the base of the expression (17) for $C_{2p-1}$ it is reduced to the
equality $[l_{j',2p}+2]+[l_{j',2p}]-[2][l_{j',2p}+1]=0$, which coincides
with the well-known identity valid for any three  successive
$q$-numbers:
$$[x+1]=[2][x]-[x-1].                             \leqno(18)
$$

\smallskip
\centerline
{\it Relations corresponding to (I.2.b)}
\smallskip

    The relation corresponding to typical vector
$\vert l_{j',2p}+2;l_{j,2p-1}+1 \rangle $ is
$$A_{2p}^{j'}(l_{j',2p}+1)A_{2p}^{j'} B_{2p-1}^j \Bigl(
\frac{B_{2p-1}^j(l_{j',2p}+2)}{B_{2p-1}^j} +
\frac{A_{2p}^{j'}(l_{j',2p}+1;l_{j,2p-1}+1)}{A_{2p}^{j'}(l_{j',2p}+1)}
\frac{A_{2p}^{j'}(l_{j,2p-1}+1)}{A_{2p}^{j'}}
$$
$$
-[2] \frac{B_{2p-1}^j(l_{j',2p}+1)}{B_{2p-1}^j}
\frac{A_{2p}^{j'}(l_{j',2p}+1;l_{j,2p-1}+1)}{A_{2p}^{j'}(l_{j',2p}+1)}
\Bigr) =0.
$$
Using explicit forms of $A_{2p}^{j'}$ and $B_{2p-1}^j$,
after simple algebra it reduces to the identity
$$
[l_{j',2p}-l_{j,2p-1}+2]+[l_{j',2p}-l_{j,2p-1}]
-[2][l_{j',2p}-l_{j,2p-1}+1]=0                      \leqno(19)
$$
which is of the same type as (18).

\smallskip
\centerline
{\it Relations corresponding to (I.3.a) }
\smallskip

The relation corresponding to typical vector
$\vert l_{j',2p}+1;l_{j'',2p}+1\rangle$ has the following form :
$$
A_{2p}^{j'}(l_{2p}^{+j''})A_{2p}^{j''}
\Bigl( 1+\frac{C_{2p-1}(l_{2p}^{+j'};l_{2p}^{+j''})}{C_{2p-1}}
-[2]\frac{C_{2p-1}(l_{2p}^{+j''})}{C_{2p-1}} \Bigr)
$$
$$                         \hskip 2mm
+A_{2p}^{j''}(l_{2p}^{+j'})A_{2p}^{j'}
\Bigl( 1+\frac{C_{2p-1}(l_{2p}^{+j'};l_{2p}^{+j''})}{C_{2p-1}}
-[2]\frac{C_{2p-1}(l_{2p}^{+j'})}{C_{2p-1}} \Bigr)
                                          = 0.           \leqno(20)
$$
Using the explicit formula (17) for $C_{2p-1}$ we transform
the expressions multiplying $A_{2p}^{j'}(l_{2p}^{+j''})A_{2p}^{j''}$ and
$A_{2p}^{j''}(l_{2p}^{+j'})A_{2p}^{j'}$ to
$-[l_{j',2p}-l_{j'',2p}-1]$ and $[l_{j',2p}-l_{j'',2p}+1]$
respectively.
Thus, the relation (20) can be rewritten as
$$
A_{2p}^{j'}A_{2p}^{j''}
\Bigl(-[l_{j',2p}-l_{j'',2p}-1]
\frac{A_{2p}^{j'}(l_{2p}^{+j''})}
{A_{2p}^{j'}}
+[l_{j',2p}-l_{j'',2p}+1]
\frac{A_{2p}^{j''}(l_{2p}^{+j'})}{A_{2p}^{j''}}
\Bigr)=0,
$$
and the latter can be verified easily with the use of explicit
expressions for $A_{2p}^{j}$ .

\medskip
\centerline
{\it Relations corresponding to (I.3.b) }
\smallskip

The proof in this case can be carried out in a  manner similar to
that of the case (I.3.a).

\smallskip
\centerline
{\it Relation corresponding to (I.1.a)}
\medskip

Next, consider the relation which corresponds to the unchanged
vector (I.1.a), i.e.,
$$
\sum_{r=1}^p \Biggl\{ (A_{2p}^r)^2 \Biggl( 2-[2]
\frac{C_{2p-1}(l_{2p}^{+r})}{C_{2p-1}}\Biggr)
                                                   \leqno(21)
$$
$$
+(A_{2p-1}^r (l_{2p-1}^{-r}))^2
\Biggl (2- [2]\frac{C_{2p-1}(l_{2p}^{-r})}{C_{2p-1}}\Biggr)
\Biggr\}=1 .
$$

Let us introduce $\phi^r$ defined as
$$
\phi^r(\{l_{1,2p+1};l_{1,2p};l_{1,2p-1}\},\ldots,
\{l_{r,2p+1};l_{r,2p};l_{r,2p-1}\},\ldots,\{l_{p,2p+1};l_{p,2p};.\})
                                                        \leqno(22)
$$
$$\equiv\frac{\prod_{s=1}^p f(l_{s,2p+1};l_{r,2p})
\prod_{s=1}^{p-1} f(l_{s,2p-1};l_{r,2p})}
{\prod_{s\ne r}^p f(l_{s,2p};l_{r,2p})f(l_{s,2p}+1;l_{r,2p})}
$$
where $f(x;y)\equiv [x+y][x-y-1]$.
\smallskip

Then, using the expression (17) for $C_{2p-1}$ and relations
$$
2-[2]\frac{[l_{r,2p}+1]}{[l_{r,2p}]}=-\frac{[2l_{r,2p}+2]}
{[l_{r,2p}+1][l_{r,2p}]} ,
$$
$$
2-[2]\frac{[l_{r,2p}-1]}{[l_{r,2p}]}=\frac{[2l_{r,2p}-2]}
{[l_{r,2p}-1][l_{r,2p}]}
$$
(the latter two can be easily verified by means of the
identity $[x+1]-[x-1]=[2x]/[x]$),
we can rewrite (21) in the form which coincides with the
relation (A.1) of Appendix if we set
$$
\phi^r(...,\{ l_{r,2p+1};l_{r,2p};l_{r,2p-1} \},...) \equiv
\frac{[2l_{r,2p}][2l_{r,2p}+2]}{[l_{r,2p}][l_{r,2p}+1]}
(A_{2p}^r(l_{r,2p}))^2 .                              \leqno(23)
$$
\bigskip
\centerline
{\it Relations corresponding to (I.1.b)}

Let us verify the following relation corresponding to typical vector
$\vert l_{j,2p-1}+1\rangle $ :
$$          \hspace{10mm}
\sum_{r=1}^p \Biggl\{
(A_{2p}^r)^2 \Biggl( 1+\frac{(A_{2p}^r (l_{2p-1}^{+j}))^2}
{(A_{2p}^r)^2}
-[2]\frac{B_{2p-1}^j (l_{2p}^{+r})}
{B_{2p-1}^j}
\frac{A_{2p}^r (l_{2p-1}^{+j})}{A_{2p}^r}\Biggr)
                                                  \leqno(24)
$$
$$
+(A_{2p}^r (l_{2p}^{-r}))^2
\Biggl ( 1+\frac{(A_{2p}^r (l_{2p}^{-r};l_{2p-1}^{+j}))^2}
{(A_{2p}^r(l_{2p}^{-r}))^2}
-[2]\frac{B_{2p-1}^j (l_{2p}^{-r})}{B_{2p-1}^j}
\frac{A_{2p}^r (l_{2p}^{-r};l_{2p-1}^{+j})}{A_{2p}^r(l_{2p}^{-r})}
\Biggr)  \Biggr\}
= 1 .
$$
Evaluating the expression at $(A_{2p}^r)^2$ we have
 $[2l_{r,2p}+2]/([l_{j,2p-1}-l_{r,2p}-1][l_{j,2p-1}+l_{r,2p}])$
(the $q$-number identities $[x][y]=[(x+y)/2]^2-[(x-y)/2]^2$
and $[x]^2-[y]^2=[x+y][x-y]$ were
utilized). As a result, the term containing $(A_{2p}^r)^2$ in
the first line of (24) takes the form
$$
\frac{1}{[2l_{r,2p}]}
[l_{r,2p}][l_{r,2p}+1]
\frac{\phi^r}{[l_{j,2p-1}-l_{r,2p}-1][l_{j,2p-1}+l_{r,2p}]}
                                                       \leqno(25)
$$
(here the correspondence (23) for $\phi^r$ and ($A^r_{2p})^2$) was
used).

Notice that in the last factor (i.e., in the ratio)
dependence on $l_{j,2p-1}$ in fact cancels out. However,
$[l_{r,2p}][l_{r,2p}+1]\equiv
-([l_{j,2p-1}-l_{r,2p}-1][l_{j,2p-1}+l_{r,2p}])\vert_{l_{j,2p-1}=0}$
and thus the expression (25) can be rewritten as
$$
\left.-\frac{1}{[2l_{r,2p}]}\phi^r\right\vert_{l_{j,2p-1}=0}.
$$
Analogously, it can be shown that the expression at
$(A_{2p}^r (l_{2p}^{-r}))^2$ in the relation (24) is equal to
$-[2l_{r,2p}-2]/([l_{j,2p-1}-l_{r,2p}][l_{j,2p-1}+l_{r,2p}-1])$
and hence the whole term containing $(A_{2p}^r (l_{2p}^{-r}))^2$
in the second line of (24) can be rewritten as
$$
\left.\frac{1}{[2l_{r,2p}]}\phi^r(l_{2p}^{-r})\right\vert_{l_{j,2p-1}=0}.
$$

Thus, the initial relation (24) is reduced to a special
case\footnote{Although the value $l_{j,2p-1}=0$ contradicts
(8) and (10), the identity (A.1) nevertheless remains true for that
value in this special case of (I.1.b).}
(with $l_{j,2p-1}=0$) of the relation (A.1)
proved for arbitrary $l_{j,2p-1}$ in the Appendix.

\smallskip
\centerline
{\it Case of (II)}
\smallskip

The proof given above (for the case of equality (I))  by
analogy carries over to the remaining cases of
equalities (II),(III),(IV).
Let us prove the equality (II).

\ni The action of RHS in (II) on the vector $\mid \xi \rangle$ leads
to the vectors :

\smallskip
(II.1) $\vert l_{j,2p}\pm 1\rangle $: \hspace{70mm} $2p$  vectors

\ni The action of LHS in (II) on the same vector leads to the vectors :

\smallskip
(II.1) $\vert l_{j,2p}\pm 1\rangle $: \hspace{70mm} $2p$  vectors

\smallskip
(II.2) $\vert l_{j',2p-1}\pm 2;l_{j,2p}\pm 1\rangle$:
                             \hspace{42mm}     $4p(p-1)$  vectors

\smallskip
(II.3) $\vert l_{j',2p-1}\pm 1;l_{j'',2p-1}\pm 1;l_{j,2p}\pm 1\rangle$ :
                            \hspace{12mm}    $4p(p-1)(p-2)$  vectors

\smallskip
(II.4) $\vert l_{j',2p-1}\pm 1;l_{j,2p}\pm 1\rangle$:
                           \hspace{44mm}$     4p(p-1)$  vectors

\smallskip
\ni(only the changed component(s) off those labelling the vectors
are indicated; the signs $\pm$ in the last three items take
their values independently).
\medskip

Let us verify the relations, which can be obtained by equating the
coefficients at the same vectors from the LHS and RHS .

\smallskip
\centerline
{\it Relations corresponding to (II.2), (II.3)}
\smallskip

Relations corresponding to (II.2) and (II.3)
can be verified in a way similar to the cases (I.2.b) and (I.3.b)
respectively.

\centerline
{\it Relations corresponding to (II.4)}
\medskip

Let us consider relation which corresponds to typical
vector $\vert l_{j',2p-1}+1;l_{j,2p}+1\rangle$.
Its explicit form is
$$
 C_{2p-1} A_{2p}^j B_{2p-1}^{j'} \Biggl (
\frac{A_{2p}^j(l_{2p-1}^{+j'})}{A_{2p}^j} +
\frac{B_{2p-1}^{j'}(l_{2p}^{+j})}{B_{2p-1}^{j'}}
\frac{C_{2p-1}(l_{2p}^{+j})}{C_{2p-1}}
- [2] \frac{B_{2p-1}^{j'}(l_{2p}^{+j})}{B_{2p-1}^{j'}}
\Biggr)
$$
$$ + C_{2p-1}(l_{2p-1}^{+j'}) A_{2p}^j B_{2p-1}^{j'} \Biggl (
\frac{A_{2p}^j(l_{2p-1}^{+j'})}{A_{2p}^j}
+ \frac{B_{2p-1}^{j'}(l_{2p}^{+j}) }{B_{2p-1}^{j'}}
\frac{C_{2p-1}(l_{2p}^{+j},l_{2p-1}^{+j'})}{C_{2p-1}(l_{2p-1}^{+j'})}
$$
$$
- [2] \frac{A_{2p}^j(l_{2p-1}^{+j'})}{A_{2p}^j}
\frac{C_{2p-1}(l_{2p}^{+j},l_{2p-1}^{+j'})}{C_{2p-1}(l_{2p-1}^{+j'})}
\Biggr) = 0.                                            \leqno(26)
$$

After evaluation (which uses explicit form of
$ A_{2p}^j, B_{2p-1}^{j'}, C_{2p-1}$ and some
q-number identities) we obtain that the expression which gives the
coefficient at $C_{2p-1} A_{2p}^j B_{2p-1}^{j'}$ is
$$
-\frac{[l_{j',2p-1}-1]}{[l_{j,2p}] ([l_{j',2p-1}-l_{j,2p}-1]
[l_{j',2p-1}-l_{j,2p}])^{1/2}}
$$
while the coefficient at $C_{2p-1}(l_{2p-1}^{+j'})
A_{2p}^j B_{2p-1}^{j'}$ is $$
\frac{[l_{j',2p-1}+1]}{[l_{j,2p}] ([l_{j',2p-1}-l_{j,2p}-1]
[l_{j',2p-1}-l_{j,2p}])^{1/2}} .
$$
Using these results, it is easy to verify validity of (26).

\medskip
\centerline
{\it Relations corresponding to (II.1)}
\smallskip

Let us consider the relation which corresponds to a
typical
vector $\vert l_{j,2p}+1\rangle$. Its explicit form is as follows :
$$
\sum_{r=1}^{p-1} \Biggl\{
(B_{2p-1}^r)^2 \Biggl( 1+\frac{(B_{2p-1}^r (l_{2p}^{+j}))^2}
{(B_{2p-1}^r)^2}
-[2]\frac{B_{2p-1}^r (l_{2p}^{+j})}
{B_{2p-1}^r}
\frac{A_{2p}^j (l_{2p-1}^{+r})}{A_{2p}^j}\Biggr) 
$$
$$    
+(B_{2p-1}^r (l_{2p-1}^{-r}))^2
\Biggl ( 1+\frac{(B_{2p-1}^r (l_{2p-1}^{-r};l_{2p}^{+j}))^2}
{(B_{2p-1}^r(l_{2p-1}^{-r}))^2}
-[2]\frac{B_{2p-1}^r (l_{2p-1}^{-r};l_{2p}^{+j})}
{B_{2p-1}^r(l_{2p-1}^{-r})}
\frac{A_{2p}^j (l_{2p-1}^{-r})}{A_{2p}^j}\Biggr)  \Biggr\}
$$
$$                  \hspace{28mm}
+(C_{2p-1})^2\left( 1+\frac{(C_{2p-1}(l_{2p}^{+j}))^2}{(C_{2p-1})^2}
-[2]\frac{C_{2p-1}(l_{2p}^{+j})}{C_{2p-1}}\right) = 1.
                                                       \leqno(27)
$$

The coefficients in (27) at
 $(B_{2p-1}^r)^2,\ $ $(B_{2p-1}^r(l_{2p-1}^{-r}))^2$
and $(C_{2p-1})^2$
 can be reduced, respectively, to the expressions
$$
\frac{[2l_{r,2p-1}+1]}
{[l_{j,2p}-l_{r,2p-1}][l_{j,2p}+l_{r,2p-1}]},\ \
-\frac{[2l_{r,2p-1}-3]}{[l_{j,2p}-l_{r,2p-1}+1][l_{j,2p}
+l_{r,2p-1}-1]},\ \ \frac{1}{[l_{j,2p}]^2}.
$$

Then the relation (27) transforms into
$$       \hspace{12mm}
\sum_{r=1}^{p-1} \frac {1}{[2l_{r,2p-1}-1]}\Biggl\{
\frac{[2l_{r,2p-1}+1][2l_{r,2p-1}-1]}
{[l_{j,2p}-l_{r,2p-1}][l_{j,2p}+l_{r,2p-1}]}
(B_{2p-1}^r)^2
          \leqno(28)
$$
$$
-\frac{[2l_{r,2p-1}-1][2l_{r,2p-1}-3]}
{[l_{j,2p}-l_{r,2p-1}+1][l_{j,2p}+l_{r,2p-1}-1]}
(B_{2p-1}^r (l_{2p-1}^{-r}))^2
\Biggr\}
+\frac{1}{[l_{j,2p}]^2} (C_{2p-1})^2 = 1.
$$
Observe that the dependence on $l_{j,2p}$ in LHS of (28) cancels
out.

Let us show that (28) is a particular case of the relation (A.1)
proved in the Appendix.
To this end, we make replacements\footnote{Although
the value $\frac32$ for $l_{j,2p+1}$ contradicts
(8) and (10), the identity (A.1) nevertheless remains true for
that value (formal for $j<p$) in this special case of (II.1).}

\vskip 2pt
$ l_{s,2p-1} \to l_{s,2p-2}+\frac{1}{2} ,
 \ \ \ \ s=1,\ldots,p-1 ;$

$ l_{s,2p} \to l_{s,2p-1}-\frac{1}{2} ,
 \ \ \ \ \ \ s=1,\ldots,p-1 ;$

$ l_{s,2p+1} \to l_{s,2p}+\frac{1}{2} ,
 \ \ \ \ \ \ s=1,\ldots,p,\ \ \ s \not= j ;$

$ l_{j,2p+1} \to \frac{3}{2} ;\qquad
\ \ \ \ \ l_{p,2p} \to -\frac{1}{2} ;$

\vskip 4pt
\ni in the LHS of (A.1).
These replacements imply (recall that $\phi^r$ is defined in (22)):
$$
\phi^r \to -\frac{[2l_{r,2p-1}+1][2l_{r,2p-1}-1]}
{[l_{j,2p}-l_{r,2p-1}][l_{j,2p}+l_{r,2p-1}]}
(B_{2p-1}^r)^2 ,      \qquad
       \frac{1}{[2l_{r,2p}]} \to  \frac{1}{[2l_{r,2p-1}-1]}
$$
for $r=1,\ldots ,p-1,$   and
$$
\phi^p \to \frac{1}{[l_{j,2p}]^2} (C_{2p-1})^2,
   \qquad  \phi^p(l_{2p}^{-p}) \to 0,
\qquad   \frac{1}{[2l_{p,2p}]} \to -1.
$$
Therefore, $\Phi(\ldots)=1$ from (A.1) goes into (28). Thus, we have
shown that the relation (28) holds simultaneously with (A.1) (which is
proved in the Appendix).

\medskip
\centerline
{\it Case of (IV) }
\smallskip

The proof of relations which appear in conjunction with (IV)
is completely analogous to that of (I). The only difference
consists in formal replacements $l_{j,2p-1} \leftrightarrow
l_{j,2p+1},$ and then $p\to p-1$. For instance, to prove the
relation corresponding to the vector
$\vert l_{j',2p-2}+2; l_{j,2p-1}+1 \rangle ,$
one uses the identity
$[l_{j',2p-2}-l_{j,2p-1}+2]+[l_{j',2p-2}-l_{j,2p-1}]
-[2][l_{j',2p-2}-l_{j,2p-1}+1]=0 ,$  which can be obtained from the
identity (19) with the replacements just pointed out.

\smallskip
\centerline
{\it Case of (III) }
\smallskip

The proof of relations corresponding to
(vectors arising in conjunction with) (III) proceeds
in complete analogy to that of {(II)} with the appropriate
replacements $l_{j,2p} \leftrightarrow l_{j,2p-2}$ being done.

\centerline
{\it Finite dimensionality, irreducibility, $*$-property}
\medskip

Concerning proof of finite dimensionality and irreducibility of
the $U_q({\rm so}_n)$ representations presented in this paper,
all things go through in complete analogy to the classical
($q=1$) case for the considered case of ${\bf m}_n$ with all
integral or all half-integral components and the restriction
$q^N\ne1.$ Finite dimensionality is verified on the base of
explicit formulas for the representation matrix elements,
with the account of the $q$-number property that, for real $x$
and $q^N\ne1$, the equality $[x]=0$ holds only if $x=0.$ Thus,
no new zeros for the matrix elements $A_{2p}^j,$  $\ B_{2p-1}^j,$
and $C_{2p-1}$ can appear besides zeros appearing in $q=1$ case
(and corresponding to limiting values in (10)).
 In the proof of irreducibility it is important to stress that
for reals $a,\ b,$ and $q^N\ne 1,$ the equality $[a]=[b]$ cannot
hold for $a\ne b.$ With this property in mind, it is enough
to trace the arguments of the proof in classical case.

The fact that the representation operators at $q=e^h$ or
$q=e^{{\rm i}h},\ h\in {\bf R},$ satisfy
the $*$-relation (7) and thus define infinitesimally unitary
(or $*$-) representations of the algebra $U_q({\rm so}_{n})$
is easily verified if one takes in account explicit matrix
elements (14)-(17).\ \ \ $\Box$

\medskip
\centerline
{\sc 4. Summary and Outlook}
\medskip

The results presented above demonstrate unambiguously applicability
of the ($q$-analogue of) Gel'fand--Tsetlin formalism for developing
representation theory of the {\it nonstandard} $q$-deformed algebras
$U_q({\rm so}_{n}),$ contrary to the opinion existing in the
literature, see [16],[17]. Within nonstandard $q$-deformation,
one is able to construct the $q$-analogues of most general finite
dimensional irreducible ${\rm so}_n$ representations (remark
that so far, only symmetric tensor representations
have been constructed [18] for the standard deformations
 $U_q({\rm B}_r)$ and $U_q({\rm D}_r)$).
Moreover, it is possible to make use of such an effective formalism in
order to construct and analyze  infinite dimensional representations
of the noncompact counterparts $U_q({\rm so}_{n-1,1})$ of
the orthogonal
$q$-algebras. The same concerns $q$-deformed inhomogeneous
(Euclidean) algebras $U_q({\rm iso}_{n})$ for which the chain
$U_q({\rm iso}_{n}) \supset U_q({\rm so}_{n}) \supset
U_q({\rm so}_{n-1})
\supset \cdots \supset U_q({\rm so}_3)$
is valid, as well as the chain
$U_q({\rm iso}_{n}) \supset U_q({\rm iso}_{n-1}) \supset
U_q({\rm iso}_{n-2})
\supset \cdots \supset U_q({\rm iso}_2).$
For the case of $U_q({\rm iso}_{n}),$ the class 1 representations
were studied
in [19], and more general representations in [20].

Analysis of representations of nonstandard $q$-deformed algebras
$U_q({\rm so}_{n})$ in the situation when the numbers characterizing
representations are not
necessarily integers or half-integers, as well as in the case of
$q$ equal to a root of unity, will be presented elsewhere.


\medskip
\centerline
{\sc Appendix}
\medskip

The aim of this appendix is to prove an important identity.

\smallskip
\ni {\bf Proposition A.} {\it The equality}
$$
\Phi(\{l_{1,2p+1};l_{1,2p};l_{1,2p-1}\},\ldots,
\{l_{p-1,2p+1};l_{p-1,2p};l_{p-1,2p-1}\},\{l_{p,2p+1};l_{p,2p};.\})
$$
$$
\equiv
\sum_{r=1}^p\frac1{[2l_{r,2p}]}
\Bigl(-\phi^r(...,\{l_{r,2p+1};l_{r,2p};l_{r,2p-1}\},...)\hspace{10mm}
\leqno(A.1)
$$
$$                     \hspace{10mm}
+
\phi^r(...,\{l_{r,2p+1};l_{r,2p}-1;l_{r,2p-1}\},...)\Bigr)
 = 1
$$
{\it with $\phi^r$ and $f(x;y)$ defined in (22) holds identically if the
$l$-coordinates (see (11) ) $l_{r,2p+1};l_{r,2p};l_{r,2p-1}$ are
consistent with inequalities (8),(10). }

\medskip
\ni{\it Proof.}
Let us show that this statement for arbitrary $l_{1,2p+1}$ is a
consequence of the statement reformulated
for $l_{1,2p+1}=l_{1,2p}+1$
( if $l_{1,2p+1}=l_{1,2p}+1$ from the very beginning, the
procedure that follows is not necessary). Using the identity
$f(x;z)-f(y;z)=
[x+z][x-z-1]-[y+z][y-z-1]=[x+y-1][x-y]=f(x;y-1)$ and its special
case $f(l_{1,2p+1};l_{r,2p})-f(l_{1,2p}+1;l_{r,2p})=
f(l_{1,2p+1};l_{1,2p})$ we get the following relations:
$$
\Phi(\{l_{1,2p+1};l_{1,2p};l_{1,2p-1}\},\ldots)-
\Phi(\{l_{1,2p}+1;l_{1,2p};l_{1,2p-1}\},\ldots)
$$
$$
\hspace{25mm}   =f(l_{1,2p+1};l_{1,2p})
\Phi(\{.;l_{1,2p};l_{1,2p-1}\},\ldots),
                                                 \leqno(A.2)
$$
$$
\Phi(\{l_{1,2p}+1;l_{1,2p};l_{1,2p-1}\},\ldots)-
\Phi(\{l_{1,2p};l_{1,2p};l_{1,2p-1}\},\ldots)
$$
$$
\hspace{16mm}   =[2l_{1,2p}] \Phi(\{.;l_{1,2p};l_{1,2p-1}\},\ldots).
                                                 \leqno(A.3)
$$

Here "." in place of $l_{1,2p+1}$ in $\Phi(\ldots)$ means that all the
multipliers depending on $l_{1,2p+1}$ are omitted.
Note that $\Phi(\{l_{1,2p};l_{1,2p};l_{1,2p-1}\},...)$
is well-defined as a formal
expression in its variables, although the value of $l_{1,2p+1}$ equal
to $l_{1,2p}$
in $\Phi(\{l_{1,2p};l_{1,2p};l_{1,2p-1}\},...)$ is excluded
by inequalities (10).

 It is possible to show by direct verification the validity of the
identity
$$
\Phi(\{l_{1,2p};l_{1,2p};l_{1,2p-1}\},\ldots)=
\Phi(\{l_{1,2p}+1;l_{1,2p};l_{1,2p-1}\},\ldots)
\vert_{l_{1,2p}\to l_{1,2p}+1}.                     \leqno(A.4)
$$
We have in LHS of (A.4)
$$
\phi^1(\{l_{1,2p};l_{1,2p};l_{1,2p-1}\},\ldots)
=\frac{\prod_{s=2}^p f(l_{s,2p+1};l_{1,2p})
\prod_{s=1}^{p-1} f(l_{s,2p-1};l_{1,2p}) f(l_{1,2p};l_{1,2p})}
{\prod_{s=2}^p f(l_{s,2p};l_{1,2p})f(l_{s,2p}+1;l_{1,2p})},
$$
\vskip 1pt
$$
\phi^1(\{l_{1,2p};l_{1,2p}-1;l_{1,2p-1}\},\ldots)=0,  \hspace{68mm}
$$
$$
\phi^r(\{l_{1,2p};l_{1,2p};l_{1,2p-1}\},\ldots
\{l_{r,2p};l_{r,2p};l_{r,2p-1}\},\ldots)             \hspace{48mm}
$$
$$
=\frac{\prod_{s=2}^p f(l_{s,2p+1};l_{r,2p})
\prod_{s=1}^{p-1} f(l_{s,2p-1};l_{r,2p}) f(l_{1,2p};l_{r,2p})}
{\prod_{ s \ge 2 , s \ne r }^p
 f(l_{s,2p};l_{r,2p})f(l_{s,2p}+1;l_{r,2p})
f(l_{1,2p};l_{r,2p})f(l_{1,2p}+1;l_{r,2p})},
$$
\vskip 1pt
\ni and in RHS of (A.4)
$$
\phi^1(\{l_{1,2p}+1;l_{1,2p};l_{1,2p-1}\},\ldots) = 0\
\stackrel {l_{1,2p}\to l_{1,2p}+1}
{\frac{\hspace{18mm}}{\hspace{18mm}}\!\!\!\longrightarrow} \ 0 ,
                                                    \hspace{31mm}
$$
\vskip 1pt
$$
\phi^1(\{l_{1,2p}+1;l_{1,2p}-1;l_{1,2p-1}\},\ldots)   \hspace{60mm}
$$
$$ \hspace{13mm}
=\frac{\prod_{s=2}^p f(l_{s,2p+1};l_{1,2p}-1)
\prod_{s=1}^{p-1} f(l_{s,2p-1};l_{1,2p}-1) f(l_{1,2p}+1;l_{1,2p}-1)}
{\prod_{s=2}^p f(l_{s,2p};l_{1,2p}-1)f(l_{s,2p}+1;l_{1,2p}-1)
}
$$
$$\hspace{15mm}
\stackrel {l_{1,2p}\to l_{1,2p}+1}
{\frac{\hspace{18mm}}{\hspace{18mm}}\!\!\!\longrightarrow}\
\frac{\prod_{s=2}^p f(l_{s,2p+1};l_{1,2p})
\prod_{s=1}^{p-1} f(l_{s,2p-1};l_{1,2p}) f(l_{1,2p}+2;l_{1,2p})}
{\prod_{s=2}^p f(l_{s,2p};l_{1,2p})f(l_{s,2p}+1;l_{1,2p})}   ,
$$
\vskip 1pt
$$
\phi^r(\{l_{1,2p}+1;l_{1,2p};l_{1,2p-1}\},\ldots,
\{l_{r,2p};l_{r,2p};l_{r,2p-1}\},\ldots)                \hspace{28mm}
$$
$$
\hspace{15mm}=\frac{\prod_{s=2}^p f(l_{s,2p+1};l_{r,2p})
\prod_{s=1}^{p-1} f(l_{s,2p-1};l_{r,2p}) f(l_{1,2p}+1;l_{r,2p})}
{\prod_{ s \ge 2 , s \ne r }^p
 f(l_{s,2p};l_{r,2p})f(l_{s,2p}+1;l_{r,2p})
f(l_{1,2p};l_{r,2p})f(l_{1,2p}+1;l_{r,2p})}
$$
$$\hspace{15mm}
\stackrel {l_{1,2p}\to l_{1,2p}+1}
{\frac{\hspace{18mm}}{\hspace{18mm}}\!\!\!\longrightarrow}\
\frac{\prod_{s=2}^p f(l_{s,2p+1};l_{r,2p})
\prod_{s=1}^{p-1} f(l_{s,2p-1};l_{r,2p}) }
{\prod_{ s \ge 2 , s \ne r }^p
 f(l_{s,2p};l_{r,2p})f(l_{s,2p}+1;l_{r,2p})f(l_{1,2p}+1;l_{r,2p})}.
$$
\vskip 2pt
\ni Using $f(l_{1,2p};l_{1,2p})=-[2l_{1,2p}]$ and
$f(l_{1,2p}+2;l_{1,2p})=[2l_{1,2p}+2]$
we can verify (A.4) by comparing the terms in
LHS and RHS of (A.4) .

Our assumption that statement (A.1) is correct for
$l_{1,2p+1}=l_{1,2p}+1$ together with (A.4) imply that LHS of
(A.3) is zero and therefore in the case of nonvanishing
$[2l_{1,2p}]$ we have the equality
$\Phi(\{.;l_{1,2p};l_{1,2p-1}\},\ldots)=0$. If we
put this result in RHS of (A.2) we get
$$
\Phi(\{l_{1,2p+1};l_{1,2p};l_{1,2p-1}\},\ldots)=
\Phi(\{l_{1,2p}+1;l_{1,2p};l_{1,2p-1}\},\ldots)=1 .
$$
Thus, it is enough to prove statement (A.1) only for the case
$l_{1,2p+1}=l_{1,2p}+1$.

In similar way we can prove that the statement (A.1)
for arbitrary $l_{1,2p-1}$ is a consequence of statement for
$l_{1,2p-1}=l_{1,2p}$ (this is superfluous if $l_{1,2p-1}=l_{1,2p}$
from the beginning).
Thus, it is enough to prove the special case of statement
(A.1) which is of the form
$$
\Phi(\{l_{1,2p}+1;l_{1,2p};l_{1,2p}\},\ldots)=1.
$$
On the other hand,
$$
\Phi(\{l_{1,2p}+1;l_{1,2p};l_{1,2p}\},
\{l_{2,2p+1};l_{2,2p};l_{2,2p-1}\},\ldots)=
\Phi(\{l_{2,2p+1};l_{2,2p};l_{2,2p-1}\},\ldots)
$$
where  all multipliers depending on $l_{1,2p+1};l_{1,2p};l_{1,2p-1}$
in the RHS are absent. This follows from the fact that
$$
\phi^1(\{l_{1,2p}+1;l_{1,2p};l_{1,2p}\},\ldots)=0,     \hspace{40mm}
$$
\vskip 1pt
$$
\phi^1(\{l_{1,2p}+1;l_{1,2p}-1;l_{1,2p}\},\ldots)=0,    \hspace{35mm}
$$
\vskip 1pt
$$
\phi^r(\{l_{1,2p}+1;l_{1,2p};l_{1,2p}\},\ldots)        \hspace{44mm}
$$
$$
=\frac{\prod_{s=2}^p f(l_{s,2p+1};l_{r,2p})
\prod_{s=1}^{p-1} f(l_{s,2p-1};l_{r,2p}) }
{\prod_{ s \ge 2 , s \ne r }^p
 f(l_{s,2p};l_{r,2p})f(l_{s,2p}+1;l_{r,2p})
}
\frac{f(l_{1,2p}+1;l_{r,2p})f(l_{1,2p};l_{r,2p})}
{f(l_{1,2p};l_{r,2p})f(l_{1,2p}+1;l_{r,2p})}
$$
in the LHS.
\vskip 2pt

Thus, fulfillment of statement (A.1) for $\ p-1\ $ triples of
variables implies
validity of statement (A.1)  for $\ p\ $ triples .
 By repeating these
arguments the analysis is
reduced to the relation with single 'truncated' triple
$\{l_{p,2p+1};l_{p,2p};.\}$
(corresponding to the $U_q({\rm so}_3)$ case) :
$$
\Phi(\{l_{p,2p+1};l_{p,2p};.\})\equiv
\frac{1}{[2l_{p,2p}]}(-f(l_{p,2p+1};l_{p,2p})
+f(l_{p,2p+1};l_{p,2p}-1))=1
$$
whose validity is verified directly .\ \ $\Box$

\bigskip
{\it Acknowledgement.} The authors are grateful to Prof. A.U.Klimyk
for discussions and valuable comments.
                The research described in this publication was made
                possible in part by CRDF Grant UP1-309 and by
                DFFD Grant 1.4/206.

\bigskip
\centerline  {\sc Bibliography}
\medskip
\small {
1. V. G. Drinfeld,{\it Hopf algebras and the quantum Yang-Baxter equation},
 Sov. Math. Dokl. {\bf 32} (1985), 254-258.

2. M. Jimbo,{\it A $q$-difference analogue of $U(g)$ and the
Yang-Baxter equation}, Lett. Math. Phys. {\bf10} (1985), 63-69.

3. I. M. Gel'fand and M. L. Tsetlin, {\it Finite dimensional
representations of the group of orthogonal matrices}, Dokl. Akad. Nauk
SSSR {\bf 71} (1950), 1017-1020.

4. I. M. Gel'fand, R. A. Minlos and Z.Ya. Shapiro,
{\it Representations of Rotation and Lorentz Groups},
Pergamon Press, New York, 1963, (Supplements).

5. L.D.Faddeev, N.Yu.Reshetikhin and L.A.Takhtadzhyan,
{ \it Quantization of Lie groups and Lie algebras},
   Leningrad Math. J. {\bf 1} (1990), 193-225.

6. A. U. Klimyk and A. M. Gavrilik, {\it Representation matrix elements
and Clebsch--Gordan coefficients of the semisimple Lie groups},
 J. Math. Phys.  {\bf 20} (1979), 1624-1642.

7. A. M. Gavrilik and A. U. Klimyk, {\it $q$-Deformed orthogonal and
   pseudoorthogonal algebras and their representations},
Lett. Math. Phys. {\bf 21} (1991), 215-220.

8. A.M.Gavrilik, I.I.Kachurik, A.U.Klimyk,
 {\it Deformed orthogonal and pseudoorthogonal Lie algebras and
their representations}, preprint ITP-90-26E, Kiev, 1990, pp.~1-17.

9. A. M. Gavrilik,
{ \it The representations of $U_q({\rm so}_{4})$ and $U_q({\rm so}_{3,1})$},
Teor. Matem. Fiz. {\bf 95} (1993), 251-257.

10. A. M. Gavrilik and A. U. Klimyk,
{ \it Representations of the $q$-deformed algebras $U_q({\rm so}_{2,1})$ and
$U_q({\rm so}_{3,1})$}, J. Math. Phys. {\bf 35} (1994), 3670-3686.

11. I. I. Kachurik and A. U. Klimyk, {\it Representations of the
$q$-deformed algebra $U'_q({\rm so}_{4})$}, J. Phys.  A {\bf 27} (1994),
 7087-7097.

12. M. Havli{\v c}ek, E. Pelantova and A. U. Klimyk,
{\it Santilli-Fairlie algebra $U_q({\rm so}_3)$: tensor products,
oscillator representations and root of unity}, Algebras, Groups and
Geometries (to appear).

13. D. B. Fairlie, {\it Quantum deformations of ${\rm SU}(2)$},
         J. Phys. A {\bf 23} (1990), L183-L187.

14. A. Odesskii, {\it An analogue of the Sklyanin algebra},
Func. Anal. Appl. {\bf 20} (1986), 152-154.

15. Yu. S. Samoilenko and L. B. Turowska, {\it Semilinear relations
and $*$-representations of deformations of ${\rm SO}(3)$},
Rep. Math. Phys. (to appear).

16. A. Chakrabarti,
{ \it ${\rm SO}(5)_q$ and Contraction: Chevalley basis representations for
$q$ generic and root of unity}, J. Math. Phys. {\bf 35} (1994),
4247-4267.

17. B. Abdesselam, D. Arnaudon and A. Chakrabarti,
{ \it Representations of $U_q({\rm SO}(5))$ and non-minimal $q$-deformation},
J.Phys.A:Math.Gen. {\bf 28} (1995), 3701-3708.

18. T. Nakashima, {\it A basis of symmetric tensor representations
for the quantum analogue of the Lie algebras ${\rm B}_n, \ {\rm C}_n$
and $\ {\rm D}_n$}, Publ. Res. Inst. Math. Sci. {\bf 26} (1990), 723-733.

19. V. A. Groza, I. I. Kachurik and A. U. Klimyk, {\it $q$-Deformed
Euclidean algebras and their representations},
Teor. Mat. Fizika {\bf 103} (1995), 467-475.

20. A. M. Gavrilik and N. Z. Iorgov, {\it $q$-Deformed inhomogeneous
algebras $U_q({\rm iso}_n)$ and their representations}, in Proc.
 of Int. Conf. "Symmetry in Nonlinear Mathematical Physics" dedicated
to memory of Wilhelm Fushchych (Kiev, July 1997).  }

\bigskip

\end{document}